\documentclass{amsart}

\usepackage[utf8]{inputenc}
\usepackage{amsmath}
\usepackage{amsfonts}
\usepackage{amssymb}
\usepackage{amsthm}
\usepackage{dsfont}
\usepackage{eufrak}
\usepackage{enumerate}

\usepackage[all]{xy}
 \usepackage[usenames,dvipsnames]{pstricks}
\usepackage{ebezier}
\usepackage[pagewise]{lineno}
\usepackage[colorlinks]{hyperref}

\usepackage{latexsym}\usepackage{amsfonts}\usepackage{amssymb}
\usepackage{newlfont}\usepackage{graphicx}\usepackage{doc}
\usepackage{color}\usepackage{multicol}\usepackage{epic}\usepackage{eepic}
\usepackage{ebezier}

\newtheorem{thrm}{Theorem}

\newtheorem{cor}[thrm]{Corollary}

\newtheorem{prop}[thrm]{Proposition}
\newtheorem{exam}[thrm]{Example}

\newcommand{\Ham}{\operatorname{Ham}}
\newcommand{\Id}{\operatorname{Id}}
\newcommand{\PSS}{\operatorname{PSS}}
\newcommand{\im}{\operatorname{im}}

\newcommand{\supp}{\operatorname{supp}}

\def \Dj{\mbox{\raise0.3ex\hbox{-}\kern-0.4em D}}

\begin{document}

\title{A note on partial quasi-morphisms and products in Lagrangian Floer homology in cotangent bundles}

\author{Jelena Kati\'c, Darko Milinkovi\'c, Jovana Nikoli\'c}

\address{Matemati\v{c}ki fakultet, Studentski trg 16, 11000 Belgrade, Serbia}

\email{jelenak@matf.bg.ac.rs, milinko@matf.bg.ac.rs, jovanadj@matf.bg.ac.rs}

\maketitle

\begin{abstract} We define partial quasi-morphisms on the group of Hamiltonian diffeomorphisms of the cotangent bundle using the spectral invariants in Lagrangian Floer homology with conormal boundary conditions, where the product compatible with the PSS isomorphism and the homological intersection product is lacking. 
\end{abstract}

\section{Introduction}\label{sec:introduction}
A {\it quasi-morphism} on a group $\mathcal{G}$  is a map
$\mu:\mathcal{G}\rightarrow\mathbb{R}$  such that there exists constant $C\geq 0$ such that
$|\mu(gh)-\mu(g)-\mu(h)|\leq C \;\; \text{for all} \;\; g,h\in\mathcal{G}$.   A quasi-morphism is {\it homogeneous} if
$\mu(g^k)=k\cdot\mu(g)$ for all $g\in\mathcal{G}$, $k\in\mathbb{Z}$.

Every quasi-morphism can be homogenized by defining
$$
\bar{\mu}(g):=\lim_{n\to\infty}\frac{\mu(g^n)}{n}.
$$
The above limit exists and $\bar{\mu}$ is a homogeneous quasi-morphism (see~\cite{PolRos}). The quasi-morphism $\bar{\mu}$ is called {\it the homogenization} of $\mu$.

If $\mathcal{G}$ is a Lie group and $\mathfrak{g}$ its Lie algebra, the derivative of a quasi-morphism is a mapping $\zeta:\mathfrak{g}\rightarrow\mathbb{R}$ called {\it a quasi-state} (see~\cite{entov,EP2,PolRos}).

If $\mathrm{Ham}(P)$ is the group of Hamiltonian diffeomorphisms of a symplectic manifold $P$,  it is known that a quasi-morphism on its universal cover $\widetilde{\mathrm{Ham}}(P)$ sometimes can be defined via {\it subadditive spectral invariants}  (see~\cite{PolRos} and the references therein).  The subadditive spectral invariant, as defined in~\cite{PolRos}, is a function $c:\widetilde{\mathrm{Ham}}(M)\rightarrow\mathbb{R}$ satisfying the axioms:
\begin{itemize}
\item[(i)] $c(\psi\phi\psi^{-1})=c(\phi)$
\item[(ii)] $c(\phi\psi)\leq c(\phi)+c(\psi)$
\item[(iii)] $\int_0^1\min(F_t-G_t)\,dt \leq  c(\phi)-c(\psi) \leq \int_0^1\max(F_t-G_t)\,dt$, where $\phi$ and $\psi$ are generated by normalized Hamiltonians $F$ and $G$
\item[(iv)] $c(\phi)\in\mathrm{Spec}(\phi)$ for nondegenerate $\phi$.
\end{itemize}
The function defined by $\nu_c(\phi):=c(\phi)+c(\phi^{-1})$ for $\phi\neq\mathrm{id}$ and $\nu_c(\mathrm{id})=0$ is a pseudo-norm on $\widetilde{\mathrm{Ham}}(M)$, called {\it the spectral pseudo-norm}.  This pseudo-norm is bounded if and only if $c$ is a quasi-morphism (Proposition 3.5.1 in~\cite{PolRos}).  If a closed symplectic manifold $M$ admits a subadditive spectral invariant $c$ with bounded spectral pseudo-norm, then $\sigma(\phi):=\lim_{n\to\infty}c(\phi^n)/n$ is a homogeneous quasi-morphism (Proposition 4.8.1 in~\cite{PolRos} ), and the normalized quasi-morphism $\mu:=-\mathrm{Vol(P)}\cdot\sigma$ also satisfies {\it the Calabi property}: for every open displaceable subset $U\subset P$ the restriction $\mu|_{\widetilde{\mathrm{Ham}}(U)}$ is equal to {\it the Calabi homomorphism}
$$
\mathrm{Cal}_U:\widetilde{\mathrm{Ham}}(U)\rightarrow\mathbb{R}, \quad \mathrm{Cal}_U(\phi):=\int_0^1\int_MH_t\,\omega^{\wedge n}\,dt.
$$

The first construction of spectral invariants is given by Viterbo~\cite{viterbo1}, by means of generating functions, and later by Oh~\cite{oh1,oh2} and Schwarz~\cite{schwarz1}. In maximal generality,  for every $\phi\in\widetilde{\mathrm{Ham}}(P)$ and for every non-zero {\it quantum} homology class $\alpha$ spectral numbers $\rho(\alpha,\phi)$ are defined as certain minimax values in Floer homology for periodic orbits (see~\cite{oh3,oh4,usher1,usher2}). They satisfy the triangle inequality 
\begin{equation}\label{eq:qtriangle}
\rho(\alpha\star\beta,\phi\psi)\leq\rho(\alpha,\phi)+\rho(\beta,\psi), 
\end{equation}
where $\star$ is the quantum product. If there is an idempotent non-zero element in the quantum homology ring, i.e. an element $e$ such that $e^2=e$, then it follows from~(\ref{eq:qtriangle}) that $c(\cdot):=\rho(e,\cdot)$ satisfies the subadditivity axiom (ii).  If a corresponding spectral pseudo-norm is bounded,  then there exists a homogeneous quasi-morphism $\mu$ on $\widetilde{\mathrm{Ham}}(M)$, which satisfies the Calabi property (see Theorem 3.1 in~\cite{entov} and Corollary 12.6.2 in~\cite{PolRos}). 
In more general case,  if $e$ is an idempotent non-zero element in a quantum homology ring,  the mapping $\mu$ constructed above has weaker properties (see Theorem 3.2 in~\cite{entov}): 
\begin{itemize}
\item instead of being quasi-additive it is {\it partially quasi-additive}, i.e. for a displaceable open set $U$ there exists $C>0$ such that 
$$
|\mu(\phi\psi)-\mu(\phi)-\mu(\psi)|\leq C\min\{\|\phi\|_U,\|\psi\|_U\},
$$
where $\|\cdot\|_U$ is {\it Banyaga's fragmentation norm} (see~\cite{banyaga}), and 
\item insted of being homogeneous it is {\it partially homogeneous}, i.e. $\mu(\phi^n)=n\mu(\phi)$ for any {\it non-negative} integer $n\in\mathbb{Z}_{\geq0}$. 
\end{itemize}
Such $\mu$ is called {\it a partial quasi-morphism}.

For a closed manifold $M$, homology class $a$ and a compactly supported Hamiltonian diffeomorpism $\phi:T^*M\rightarrow T^*M$,  the spectral numbers $\ell (a,\phi)$ can be defined via Lagrangian Floer homology~\cite{oh1,oh2}. They satisfy the triangle inequality  
\begin{equation}\label{eq:ltriangle}
\ell(\alpha\cap\beta,\phi\psi)\leq \ell(\alpha,\phi)+\ell(\beta,\psi), 
\end{equation}
where $\cap$ is the intersection product in homology.  Monzner, Vichery and Zapolsky~\cite{MVZ} defined the mapping $\mu_0:\mathrm{Ham}(T^*M)\rightarrow\mathbb{R}$ on the group of compactly supported Hamiltonian diffeomorphisms in $T^*M$ by\footnote{More generally, they defined the mapping $\mu_a$ for each $a\in H^1(M)$, and proved that it has the properties anologous to those of a partial quasi-morphism (due to the later results of Shelukhin~\cite{shelukhin1} and Kislev--Shelukhin\cite{KS}, $\mu_a$ give rise to genuine quasi-morphisms for some $M$). If $M=\mathbb{T}^n$,   $\mu_p(\phi_1^H)=\overline{H}(p)$, where $\overline{H}$ is the Viterbo's homogenization~\cite{viterbo2} of $H$  and $p\in\mathbb{R}^n\cong H^1(\mathbb{T}^n)$~\cite{MVZ}  (see also~\cite{MZ}). }
$$
\mu_0(\phi):=\lim_{n\to\infty}\frac{\ell([M],\phi^n)}{n}.
$$
Since $[M]$ is an idempotent non-zero element in homology $H_*(M)$, the subadditivity here follows from~(\ref{eq:ltriangle}).

In this paper we ponder over the construction of partial quasi-morphisms and quasi states by using the spectral numbers from Lagrangian Floer homology with conormal boundary conditions in cotangent bundles.   Due to the lack of product satisfying the triangle inequalities analogous to~(\ref{eq:qtriangle}) and~(\ref{eq:ltriangle}),  these spectral numbers do not satisfy the subadditivity axiom (ii). We circumvent this shortcoming by using the exterior intersection product compactible with the $\mathrm{PSS}$ isomorphism as constructed in~\cite{D}, to prove the following theorem (see Section~\ref{sec:results} for more precise formulation):

\begin{thrm}
Let $M$ be a closed  manifold and $N\subset M$ a closed submanifold.  There are well defined partial quasi-morphism
$$
\sigma^N:\mathrm{Ham}(T^*M)\rightarrow\mathbb{R},  \quad \sigma^N(\phi):=\overline{\lim}\frac{\ell^N_+(\phi^n)}{n}
$$
and partial quasi-state
$$
\zeta^N:C^{\infty}_c(T^*M)\rightarrow\mathbb{R}, \quad \zeta^N(H):=\sigma^N(\phi^H_1),
$$
where $\ell^N_+(\cdot):=\ell^N_+([N],\cdot)$ is a spectral number defined as a minimax value in Lagrangian Floer homology $HF_*(o_M,\nu^*N:H)$.
\end{thrm}

In Section~\ref{sec:results} we give necessary definitions and precise formulation of the previous theorem.  In Section~\ref{sec:discussion} we discuss the existence of the limit in previous theorem and the products in Lagrangian Floer homology with conormal boundary conditions, and we show that our partial quasi-morphisms are non-trivial and different from those in~\cite{MVZ}.  In Section~\ref{sec:proof} we prove the results from Section~\ref{sec:results}. In Section~\ref{sec:boundary}, we briefly discuss a version of the previous theorem for the case where $N$ is a manifold with boundary.

\bigskip 

\textbf{Acknowledgement.} The third author is grateful to Octav Cornea for useful discussion at the workshop {\it Current trends in symplectic topology} in Montreal 2019,  and to Morimichi Kawasaki for the email discussion that helped us clarify some concepts related to this paper. The authors' research is partially supported by the contract of the Ministry of Education, Science and Technological Development of the Republic of Serbia No. 451-03-9/2021-14/200104.

\section{Preliminaries and results}\label{sec:results}
Let $M$ denote a smooth closed manifold and $N\subset M$ its closed submanifold. The cotangent bundle of $M$, $T^*M$, is endowed with a symplectic structure $\omega=-d\lambda$, where $\lambda$ is the tautological Liouville $1$-form. 

Suppose that the Hamiltonian $H:T^*M\times[0,1]\to\mathbb{R}$ is a smooth compactly supported map such that the intersection $\nu^*N\cap\phi_H^1(o_M)$ is transverse.  We denote the Floer homology for the pair $(o_M,\nu^*N)$ by $HF_*(o_M,\nu^*N:H)$.  The filtration of homology by the action functional
$$
{\mathcal A}_H(\gamma)=-\int\gamma^*\lambda+\int_0^1H(\gamma(t),t)\,dt.
$$
defines filtered Floer homology $HF_*^\lambda(o_M,\nu^*N:H)$. It is well known that Floer homology is isomorphic to the singular homology of $N$ and thus to Morse homology of $N$ (see~\cite{P}).  We denote the isomorphism of PSS type between these homologies by
\begin{equation}\label{eq:pss}
\PSS:HM_*(N)\to HF_*(o_M,\nu^*N:H)
\end{equation}
(see~\cite{D} for details). For $\alpha\in HM_*\setminus\{0\}$ we can define a conormal spectral number
$$\ell(\alpha;o_M,\nu^*N:H)=\inf\{\lambda\,|\,\PSS(\alpha)\in\im(\imath^\lambda_*)\},$$
where $$\imath^\lambda_*:HF_*^\lambda(o_M,\nu^*N:H)\to HF_*(o_M,\nu^*N:H)$$ is the inclusion induced by the inclusion map of chain complexes $$\imath^\lambda:CF_*^\lambda(H)\to CF_*(H).$$

Following~\cite{Aur} we can define a natural homology action homomorphism of $HF_*(o_M,o_M)$ on $HF_*(o_M,\nu^*N)$. Note that $HF_*(o_M,o_M)$ stands for Floer homology for conormal bundle in the special case when $M=N$. This is the standard product in Lagrangian Floer homology. Moreover, we can relate it, via the PSS isomorphism, to the action on the Morse side where it becomes the action of $HM_*(M)$ on $HM_*(N)$ via the external intersection product. As a result we obtain a modified version of a triangle inequality for spectral invariants (see Theorem 1.3 in~\cite{D})
\begin{equation}\label{triang1pre}
\ell(\alpha\widetilde\cap\beta;o_M,\nu^*N:H_1\sharp H_2)\le \ell(\alpha;o_M,o_M:H_1)+\ell(\beta;o_M,\nu^*N:H_2).
\end{equation}
Here, $\alpha\in HM_*(M)$ and $\beta\in HM_*(N)$ are such that $\alpha\widetilde\cap\beta\neq0$ and $\widetilde\cap$ is the exterior intersection product on Morse homology (see~\cite{D}):
\begin{equation}\label{eq:extcap}
\widetilde\cap:HM_r(M)\otimes HM_s(N)\to HM_{r+s-\mathrm{dim}\,M}(N).
\end{equation}
Taking into account the dimensions of the moduli spaces involved in the definition of this product~\cite{D}, one can show that for $r=\mathrm{dim}\,M$ and the fundamental class $[M]\in HM_{\dim\,M}(M)$
$$
[M]\,\widetilde{\cap}\,\cdot: HM_s(N) \to HM_s(N).
$$
is the identity map. Therefore, we obtain the inequality
\begin{equation}\label{eq:spectriangle}
\ell([N];o_M,\nu^*N:H_1\sharp H_2)\le \ell([M];o_M,o_M:H_1)+\ell([N];o_M,\nu^*N:H_2),
\end{equation}
where $[N]$ is the fundamental class of $N$ and $[M]$ is the fundamental class of $M$.

Let us denote
$$\ell_+^N(H)=\ell([N];o_M,\nu^*N:H).$$ 
Following~\cite{MVZ} (Proposition 2.6 and Remark 2.7) we can conclude that $\ell_+^N(H)=\ell_+^N(K)$ if $\phi_H^1=\phi_K^1$.  This means that we can see these invariants as a map with a domain $\Ham(T^*M)$.  From  inequality~(\ref{eq:spectriangle}) we get the following form of a triangle inequality
\begin{equation}\label{triang1}
\ell^N_+(\phi\psi)\le \ell^M_+(\phi)+\ell^N_+(\psi),
\end{equation}
for every $\phi,\psi\in\Ham(T^*M)$.

The homogenization of spectral invariant, that corresponds to the fundemental class, is a well known way to obtain a partial quasi-morphism on the Hamiltonian group. Following~\cite{MVZ}, first we try to homogenize $\ell_+^N$. But a limit that defines homogenization
\begin{equation}\label{Nema}
\lim\limits_{n\to+\infty}\frac{\ell_+^N(\phi^n)}{n}
\end{equation}
might not exist (see Section~\ref{sec:discussion} and Example~\ref{Primer} below). It turns out that the upper limit  of this sequence exists and it satisfies well known properties of partial quasi-morphism. More precisely, if we define a map
\begin{equation}\label{ParcKvazi}
\sigma^N:\Ham(T^*M)\to\mathbb{R},\quad
\sigma^N(\phi)=\overline{\lim}\frac{\ell_+^N(\phi^n)}{n},
\end{equation}
then the following theorem holds:

\begin{thrm}\label{ThrmQuasiMorph}
The map $\sigma^N$ is well defined and it satisfies properties of partial quasi-morphisms
\begin{enumerate}[(i)]
\item\label{Homog} (Homogeneity) $\sigma^N(\phi^l)=l\cdot\sigma^N(\phi)$, for $l\in\mathbb{Z}_{\geq0}$;
\item\label{ConjInvar} (Conjugation invariance) $\sigma^N(\phi)=\sigma^N(\psi\phi\psi^{-1})$ for any $\phi,\psi\in\Ham(T^*M)$;
\item\label{Lipsic} If the Hamiltonian $H$ generates diffeomorphism $\phi$ and $K$ generates $\psi$ then it holds
$$\int_0^1\min(H_t-K_t)\,dt\le\sigma^N(\phi)-\sigma^N(\psi)\le\int_0^1\max(H_t-K_t)\,dt;
$$
\item\label{Razdvojivo} $\sigma^N(\phi)=0$ for any $\phi\in\Ham_U(T^*M)$ where $U\subset T^*M$ is a displaceable subset and $\Ham_U(T^*M)$ the group of Hamiltonian diffeomorphsms supported in $U$;
\item\label{Fragmen} If $\mathcal{U}$ is a collection of open subsets of $T^*M$ such that the spectral displacement energy $e(\cdot)$ satisfies
$$
e(\mathcal{U}):=\sup_{U\in\mathcal{U}}e(U)<+\infty
$$
then
$$
|\sigma^N(\phi\psi)-\sigma^N(\psi)|\le e(\mathcal{U})\|\phi\|_{\mathcal{U}},
$$
for every $\phi,\psi\in\Ham(T^*M)$;
\item\label{Commute} For commuting Hamiltonian diffeomorphisms $\phi$ and $\psi$ it holds $\sigma^N(\phi\psi)\le\sigma^M(\phi)+\sigma^N(\psi)$ and $\sigma^N(\phi\psi)\le\sigma^M(\psi)+\sigma^N(\phi)$.
\end{enumerate}
\end{thrm}

Note that the property~(\ref{Fragmen}) gives us the property of controlled quasi-additivity of $\sigma^N$. Indeed, if $U$ is a displaceable open subset of $T^*M$ then 
$$|\sigma^N(\phi)|\le K\|\phi\|_{U}
$$
and 
$$
|\sigma^N(\phi\psi)-\sigma^N(\psi)-\sigma^N(\phi)|\le 2K\|\phi\|_{U},
$$
where $K=e(U)$ is a positive constant depending only on $U$ .
Since for commuting $\phi$ and $\psi$ it obviously holds $\sigma^N(\phi\psi)=\sigma^N(\psi\phi)$, we also get the inequality
$$|\sigma^N(\psi\phi)-\sigma^N(\phi)-\sigma^N(\psi)|\le 2K\|\psi\|_{U}.
$$
Therefore,  we conclude
$$|\sigma^N(\phi\psi)-\sigma^N(\psi)-\sigma^N(\phi)|\le 2K\min\{\|\phi\|_{U},\|\psi\|_{U}\}.
$$

The partial quasi-morphism $\sigma^N$ defines a map $\zeta^N:C^\infty_c(T^*M)\to\mathbb{R}$
$$\zeta^N(H)=\sigma^N(\phi_H^1).
$$
It turns out that $\zeta^N$ satisfies the properties of partial symplectic quasi-states.
\begin{thrm}\label{ThrmQuasiState} The map $\zeta^N$ satisfies the properties listed below:
\begin{enumerate}[(i)]
\item\label{NormState}(Normalization) $\zeta^N(0)=0$;
\item\label{StabilState}(Stability) For every $H,K\in C^\infty_c(T^*M)$ it holds 
$$\min\limits_{T^*M}(H-K)\le\zeta^N(H)-\zeta^N(K)\le\max\limits_{T^*M}(H-K);$$
\item\label{MonotState}(Monotonicity) $\zeta^N(H)\le\zeta^N(K)$ for $H\le K$; 
\item\label{HomState}(Homogeneity) $\zeta^N(sH)=s\zeta^N(H)$ for every $s\ge0$;
\item\label{InvState}(Invariance) $\zeta^N$ is invariant under the action of $\Ham(T^*M)$ on $C^\infty_c(T^*M)$;
\item\label{VanState}(Vanishing) If the support of $H$ is a displaceable set then $\zeta^N(H)=0$; 
\item\label{QAddState}(Partial quasi-additivity) If $H,K$ are two functions such that $\{H,K\}=0$ and the support of $K$ is displaceable then
$$\zeta^N(H+K)=\zeta^N(H).
$$
\end{enumerate}
\end{thrm}

\section{Discussion}\label{sec:discussion}

\subsection{Products} As we mentioned in Section~\ref{sec:introduction},  subadditivity of spectral invariants can be derived from the triangle inequality for Floer theoretical spectral invariants and the existence of an idempotent non-zero element with respect to a multiplicative structure in the homology,  provided that such a multiplicative structure exists.  The triangle inequality for Lagrangian spectral invariants is a consequence of the behavior of the restriction of the action functional to the pairs of pants that connect the elements participating in the product defined by it~\cite{oh1,oh2}. In the case $M=N$ the pair of pants product in Lagrangian Floer homology intertwines with the intersection product in the homology of $M$ and the PSS isomorphism~(\ref{eq:pss}). The fact that the fundamental class $[M]$ is an idempotent element in relation to the  intersection product was used by Monzer, Vichery and Zapolsky~\cite{MVZ} to prove subadditivity. If $\mathrm{dim}\,N<\mathrm{dim}\,M$ this is not possible (see Corollary~\ref{cor:proizvod} below).

\begin{figure}[htb]
\includegraphics[width=6cm,height=2cm]{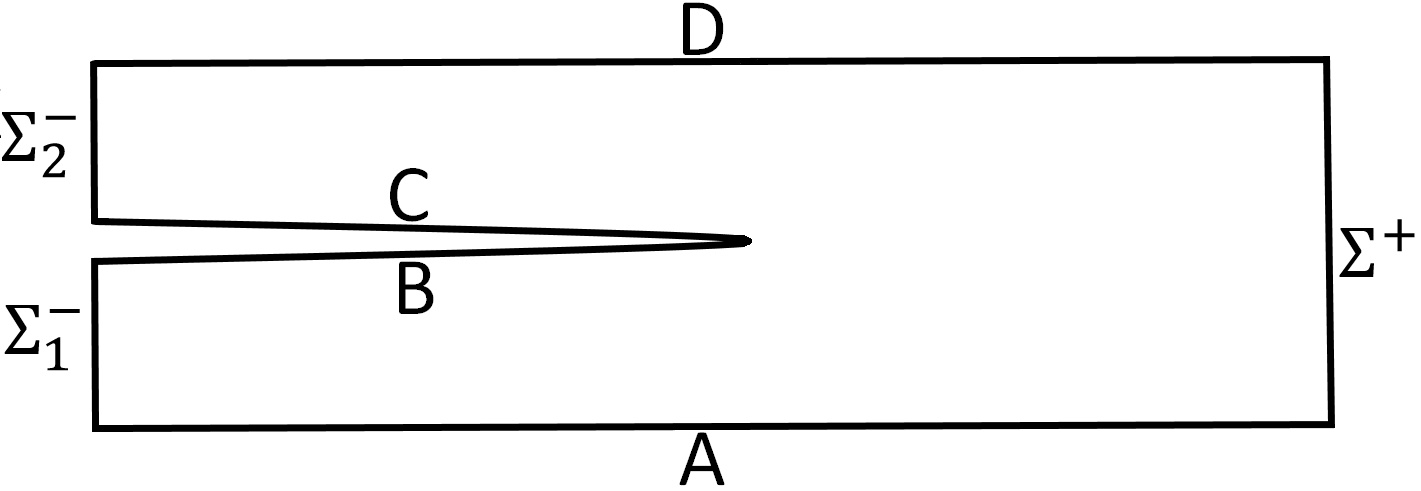}
\caption{Riemann surface $\Sigma$}
\label{Pants}
\end{figure}

Recall that the pair of pants product on $HF_*(o_M,\nu^*N:H)$ is defined by counting perturbed holomorphic pairs--of--pants with appropriate boundary conditions.  Let $\Sigma$ be a Riemann surface with boundary that has 3 strip--like ends (see Figure~\ref{Pants}): $\Sigma^-_1,\Sigma^-_2\approx [0,1]\times(-\infty,0]$, $\Sigma^+\approx [0,1]\times[0,+\infty)$. Let us denote by $\Sigma^0=\Sigma\setminus\left(\Sigma^-_1\cup\Sigma^-_2\cup\Sigma^+\right)$ the compact part of the surface $\Sigma$.

For $x_i^-\in CF_*(H_i^-)$, $i\in\{1,2\}$ and $x^+\in CF_*(H^+)$ the manifold ${\mathcal M}(x_1^-,x_2^-;x^+)$ used in this construction is the set of holomorphic maps $u$ such that
\begin{itemize}
  \item $u:\Sigma\to P$;
  \item smooth perturbation in Cauchy–-Riemann equation on the strip--like end $\Sigma_i^-$ (respectively, $\Sigma^+$) is defined by the Hamiltonian $H_i^-$ (respectively, $H^+$);
  \item there is no perturbation on the compact part of the pants $\Sigma^0$;
  \item $u$ maps lines $A$ and $C$ to $o_M$;
  \item $u$ maps lines $B$ and $D$ to $\nu^*N$ and
  \item at the strip--like end, $u$ converges (at the $\infty$) to the appropriate Hamiltonian path ($x_1^-$, $x_2^-$ or $x^+$).
\end{itemize}

More precisely,  the Riemann surface $\Sigma$ can be defined as a disjoint union ${\mathbb R}\times[-1,0]\sqcup{\mathbb R}\times[0,1]$ with identification $(s,0^-)\sim(s,0^+)$ for $s\geq0$. The above mentioned strip-like ends $\Sigma_i^-$ and $\Sigma^+$ are defined via holomorphic embeddings
\begin{equation*}
\psi_i^-:(-\infty,0]\times[0,1]\to\Sigma,\,i\in\{1,2\};\;\psi^+:[0,+\infty)\times[0,1]\to\Sigma,
\end{equation*}
as $\Sigma^-_i=\verb"Im"\,\psi_i^-$, $i\in\{1,2\}$ and $\Sigma^+=\verb"Im"\,\psi^+$. Let $\rho_R:(-\infty,0]\to[0,1]$ be a smooth function such that $\rho_R(s)=1$ for $s\le -R-1$, and $\rho_R(s)=0$, for $s\ge -R$ ($R$ is a positive fixed number). A function
$\tilde{\rho}_R:[0,+\infty)\to[0,1]$ is defined by $\tilde{\rho}_R(s)=\rho_R(-s)$.
The precise definition of the set of holomorphic pair-of-pants is
\begin{equation*}\label{eq:moduli}
{\mathcal M}(x_1^-,x_2^-;x^+)=\left\{
\begin{array}{ll}
u:\Sigma\to P, \\
\frac{\partial u_i^-}{\partial s}+J\frac{\partial u_i^-}{\partial t}=\rho_RJX_{H_i^-}(u_i^-),\,u_i^-=u\circ\psi_i^-,\,i\in\{1,2\},\\
\frac{\partial u^+}{\partial s}+J\frac{\partial u^+}{\partial t}=\tilde{\rho}_RJX_{H^+}(u^+),\,u^+=u\circ\psi^+,\\
\frac{\partial u^o}{\partial s}+J\frac{\partial u^o}{\partial t}= 0,\,u^o=u|_{\Sigma^o},\\
u(s,-1)\in o_M, u(s,1)\in \nu^*N,\,s\in{\mathbb{R}},\\
u(s,0^-)\in \nu^*N, u(s,0^+)\in o_M,\, s\leq0,\\
u_i^-(-\infty,t)=x_i^-(t),\,i\in\{1,2\};\;u^+(+\infty,t)=x^+(t).
\end{array}
\right.
\end{equation*}
As we can see, on the slit of the pants we have a jump. One part of the boundary (line $B$ on the Figure \ref{Pants}) goes to $\nu^*N$ and the other part (line $C$ on the Figure \ref{Pants}) goes to $o_M$. Holomorphic strips with jumping boundary conditions appear in \cite{AS}.   From Fredholm analysis it follows that ${\mathcal M}(x_1^-,x_2^-;x^+)$ has a structure of a manifold with corners (see Section 5 in \cite{D} for more details). 

Our product is defined on generators of $CF_*$ by
\begin{equation*}
  x_1^-\ast x_2^-=\sum_{x^+}\sharp_2{\mathcal M}(x_1^-,x_2^-;x^+)\,x^+.
\end{equation*}
Here, $\sharp_2{\mathcal M}(x_1^-,x_2^-;x^+)$ denotes the (mod 2) number of elements of a zero dimensional component of ${\mathcal M}(x_1^-,x_2^-;x^+)$. We extend the product $\ast$ by bilinearity and obtain the product on the homology
\begin{equation*}
  \ast:HF_*(o_M,\nu^*N:H_1^-)\otimes HF_*(o_M,\nu^*N:H_2^-)\to HF_*(o_M,\nu^*N:H^+).
\end{equation*}

\begin{prop}
For generic compactly supported Hamiltonians $H_1^-$,$H_2^-$ and $H^+$ the moduli space ${\mathcal M}(x_1^-,x_2^-;x^+)$ is a manifold with corners of dimension
\begin{equation*}
  \dim{\mathcal M}(x_1^-,x_2^-;x^+)=\mu_N(x_1^-)+\mu_N(x_2^-)-\mu_N(x^+)+\frac{1}{2}\dim N-\dim M.
\end{equation*}
\end{prop}

\begin{figure}[htb]
\includegraphics[width=11cm,height=3cm]{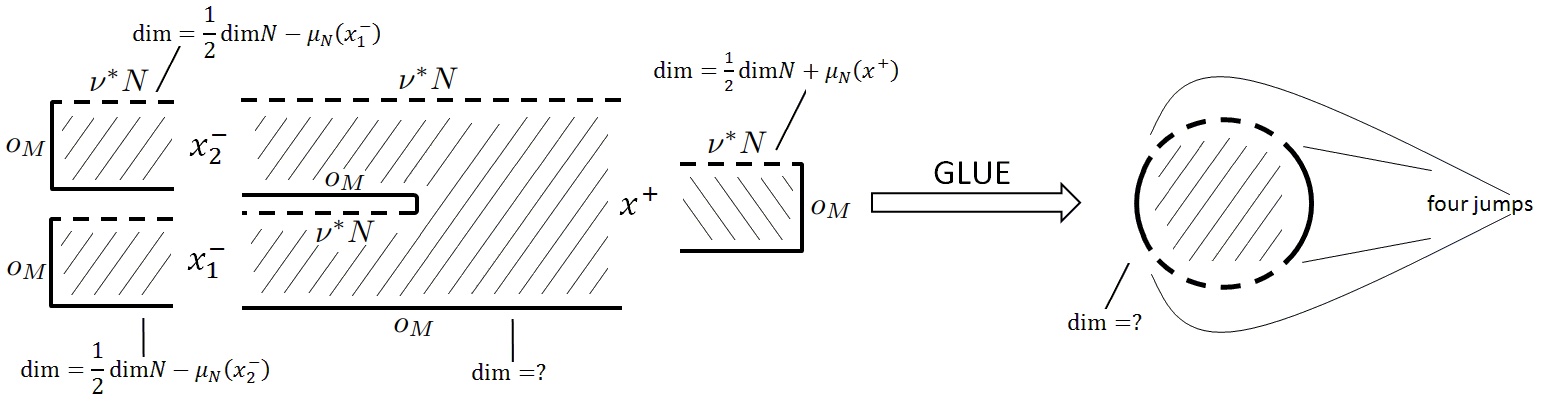}
\caption{Gluing of half strips and a pair-of-pants}
\label{Glue1}
\end{figure}

\begin{figure}[htb]
\includegraphics[width=9cm,height=2cm]{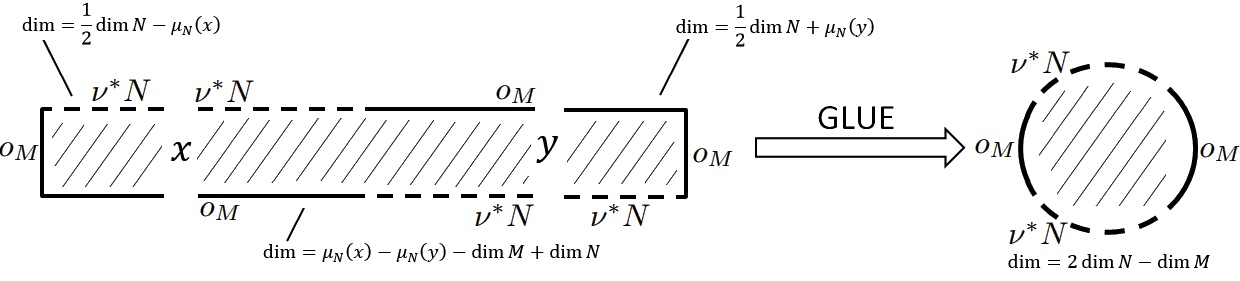}
\caption{Gluing of half strips and a whole strip}
\label{Glue2}
\end{figure}

{\it Proof:} The dimension of $\mathcal{M}(x_1^-,x_2^-;x^+)$ is equal to the Fredholm index of the linearized Cauchy--Riemann operator at $u\in\mathcal{M}(x_1^-,x_2^-;x^+)$. The easiest way to compute this index is by using the gluing of Fredholm operators. We cap the half strips whose asymptotic limits are $x_i^-,x^+$ at the appropriate strip--like ends of $\Sigma$ (see Figure \ref{Glue1}). That way we obtain some new Fredholm operator. We can compute the index of the new operator if we cap two half strips to the whole strip that has one jump at $\mathbb R\times\{0\}$ and one jump at $\mathbb R\times\{1\}$ (see Figure \ref{Glue2}). A similar idea was given in \cite{oh2} and \cite{schwarz2} in a different context.

The moduli space of half strips is defined as
\begin{equation*}
W^s(x,H)=\left\{
\begin{array}{ll}
u:[0,+\infty)\times[0,1]\rightarrow P,\;\frac{\partial u}{\partial s}+J\frac{\partial u}{\partial t}=\tilde{\rho}_RJX_{H}(u), \\
u(s,0)\in o_M, u(s,1)\in \nu^*N, u(0,t)\in o_M, s\geq 0, t\in[0,1], \\
u(+\infty,t)=x(t).
\end{array}
\right.
\end{equation*}
The dimension of $W^s(x,H)$ is computed in~\cite{AS}, $\dim W^s(x,H)=\frac{1}{2}\dim N-\mu_N(x)$.
The dimension of the space of the opposite half strips, $(-\infty,0]\times[0,1]$, is $\frac{1}{2}\dim N+\mu_N(x)$.

The moduli space of whole strips is defined by
\begin{equation*}
\overline{{\mathcal M}}(x,y)=\left\{
\begin{array}{ll}
u:{\mathbb R}\times [0,1]\rightarrow P,\;\frac{\partial u}{\partial s}+J\frac{\partial u}{\partial t}= \beta(s)JX_H(u),\\
u(s, 0),u(-s, 1)\in o_M, \; u(s, 1),u(-s, 0)\in \nu^*N,\;s\leq0,\\
u(-\infty,t)=x(t), \;
u(+\infty,t)=y(t).
\end{array}
\right.
\end{equation*}
Here, $\beta(s)$ is a smooth function such that $\beta(s)=1$ for $|s|\geq R+1$ and $\beta(s)=0$ for $|s|\leq R$.  The dimension of $\overline{{\mathcal M}}(x,y)$ is computed in \cite{AS}; $\dim \overline{{\mathcal M}}(x,y)=\mu_N(x)-\mu_N(y)-\dim M+\dim N$.
Now the proof follows from the relation
\begin{equation*}
\begin{aligned}
  &\left(\frac{1}{2}\dim N-\mu_N(x_1^-)\right)+\left(\frac{1}{2}\dim N-\mu_N(x_2^-)\right)+\dim{\mathcal M}+\left(\frac{1}{2}\dim N+\mu_N(x^+)\right)= \\
  &=\left(\frac{1}{2}\dim N-\mu_N(x)\right)
  +\left(\mu_N(x)-\mu_N(y)-\dim M+\dim N\right)+\left(\frac{1}{2}\dim N+\mu_N(y)\right).
\end{aligned}
\end{equation*}
\qed

\begin{cor}\label{cor:proizvod}
The degree of the product $\ast$ is $-\dim M$, i.e.
\begin{equation*}
  \ast:HF_r(H_1^-)\otimes HF_s(H_2^-)\to HF_{r+s-\dim M}(H^+).
\end{equation*}
In particular, the homological intersection product in $H_*(N)$,
$$
\cap: H_r(N)\otimes H_s(N)\rightarrow H_{r+s-\dim\,N}(N)
$$
does not intertwine with the PSS isomorphism~(\ref{eq:pss}) and the pair of pants product $\ast$, unless $\dim\,N=\dim\,M$.
\end{cor}

\subsection{Existence of the limit}
 We have already mentioned that the limit $\lim\limits_{n\to+\infty}\frac{\ell_+^N(\phi^n)}{n}$ might not exist. Let us discuss this point briefly. Using the same notation and the idea as in~\cite{PolRos} we define the sequences
$$
a_n=n\,\ell_+^M(\phi^{-1})+\ell_+^N(\phi^n)
$$
and
$$
b_n=n\,\ell_+^M(\phi^{-1})+\ell_+^M(\phi^n).
$$
If $\lim\limits_{n\to+\infty}\frac{a_n}{n}$ exists then the limit~(\ref{Nema}) also exists, and vice versa. From the triangle inequality~(\ref{triang1}) we can conclude that it holds
\begin{equation}\label{Property1}
a_{m+n}\le a_n+b_m,
\end{equation}
for every $m,n\in\mathbb{N}$. This means that the sequence $a_n$ is not subadditive but satisfies this weaker inequality (that is not good enough for the existence of the standard homogenization). It is known that $b_n$ is subadditive sequence and from Fekete lemma it follows that $\lim\limits_{n\to+\infty}\frac{b_n}{n}$ exists. Subadditivity of $b_n$ was proven in~\cite{MVZ}. It also follows from~(\ref{triang1}) in case when $N=M$.
The same inequality~(\ref{triang1}) and the fact that the spectral invariant of the identity is zero imply the non-negativity of the sequences
\begin{equation}\label{Property2}
a_n,b_n\ge0.
\end{equation}
If the inputs  in~(\ref{triang1pre}) are $\alpha=[M]$  and $H_2=0$ we conclude that the conormal spectral invariants are bounded for every non-zero homology class $\beta\in HM_*(N)$:
\begin{equation}\label{ineq2}
\ell(\beta;o_M,\nu^*N:H)\le \ell([M];o_M,o_M:H).
\end{equation}
It follows that
\begin{equation}\label{Property3}
a_n\le b_n,\,\,n\in\mathbb{N}.
\end{equation}
We can control the difference between two elements of the sequence $a_n$. The sequence $a_n$ is non-decreasing because
\begin{equation}\label{Property4}
\begin{aligned}
a_{n+1}-a_n&=(n+1)\,\ell_+^M(\phi^{-1})+\ell_+^N(\phi^{n+1})-n\,\ell_+^M(\phi^{-1})-\ell_+^N(\phi^{n})\\
&=\ell_+^M(\phi^{-1})+\ell_+^N(\phi^{n+1})-\ell_+^N(\phi^{n})\ge0,
\end{aligned}
\end{equation}
(the last inequality follows again from~(\ref{triang1})). We want to 
estimate $a_{n+1}-a_n$ from above. 
The relation 2.  inTheorem 2.20 in~\cite{MVZ} states 
\begin{equation}\label{LipsicInvar}
\int_0^1\min(H_t-K_t)\,dt\le 
\ell(\alpha;o_M,\nu^*N:H)-\ell(\alpha;o_M,\nu^*N:K)
\le\int_0^1\max(H_t-K_t)\,dt.
\end{equation}
We also know that if a Hamiltonian $F$ generates $\phi=\phi_F^1$ and $K$ generates $\phi^n=(\phi_F^1)^n$ ($K$ is actually $F\sharp...\sharp F$, $n-$times) then 
$$
H(x,t)=K(x,t)+F(((\phi_H^t)^n)^{-1}x,t)
$$
generates $(\phi_F^t)^{n+1}=(\phi_F^t)^n\circ\phi_F^t$. If we put this in the  inequality~(\ref{LipsicInvar}) with $\alpha=[N]$  we get
$|\ell^N_+(\phi^{n+1})-\ell^N_+(\phi^n)|\le\|H\|_{\infty}$. Finally 
\begin{equation}\label{Property5}
a_{n+1}-a_n = |a_{n+1}-a_n|
\le|\ell^M_+(\phi^{-1})|+|\ell^N_+(\phi^{n+1})-\ell^N_+(\phi^n)|\le C,
\end{equation}
where $C$ is some positive constant. 

All these properties~(\ref{Property1})-(\ref{Property5}) of sequence $a_n$ are not sufficient for the existence of the limit $\lim\limits_{n\to+\infty}\frac{a_n}{n}$.

\begin{exam}\label{Primer} 
There exists a sequence $a_n$ that satisfies the properties~(\ref{Property1}) 
to~(\ref{Property5})
and yet the limit $\lim\limits_{n\to\infty}a_n/n$ does not exist.
Indeed, take $b_n:=1$ for all $n$, $C:=1$  and define
\begin{itemize}
\item  $a_1:=1$
\item $a_{n+1}:=a_n$ for all $n\in\{1,2,\ldots,n_1-1\}$, where $n_1$ is
the least positive integer $n$ such that  it holds $a_n/n<1/3$
\item $a_{n+1}:=a_n+1$ for all $n\in\{n_1,n_1+1,\ldots,n_2-1\}$,  where
$n_2\in\{n_1+1,n_1+2,\ldots\}$ is the least positive integer $n$
 satisfying $a_n/n>1/2$
 \item $a_{n+1}:=a_n$ for all $n\in\{n_2,n_2+1,\ldots,n_3-1\}$,  where
$n_3\in\{n_2+1,n_2+2,\ldots\}$ is the least positive integer $n$
 satisfying $a_n/n<1/3$
\end{itemize}
etc. It is easy to see that these two procedures (keeping $a_{n+1}:=a_n$
and setting $a_{n+1}:=a_n+1$) will alternate infinitely many time.
Thus the sequence $a_n/n$ has at least two accumulation points, since it possesses a
subsequence $a_{n_k}/n_k>1/2$ and a subsequence $a_{n_l}/n_l<1/3$.
\end{exam}

\subsection{Nontriviality and the comparison with the case $N=M$} In the following example we compute $\sigma^N$ for some Hamiltonian diffeomorphisms in order to show that our partial quasi-morphisms are not trivial.
We also point out the difference between these partial quasi-morphisms and those constructed in~\cite{MVZ}.
 
Let $M$ be $\mathbb{S}^1$ and  $f:\mathbb{S}^1\to\mathbb{R}$ a Morse function with two critical points: the minimum point $x_1$  and the maximum point $x_2$.  Let $N$ be a point $x_1$, and thus 
$\nu^*N=T^*_{x_1}\mathbb{S}^1$. 
We lift $f$ to the cotangent bundle, $H_f=f\circ\pi$ where $\pi:T^*\mathbb{S}^1\to\mathbb{S}^1$ is the projection. The flow of the (non-compactly supported) Hamiltonian $H_f$  is given by
\begin{equation}\label{eq:Ham1}
\phi_{H_f}^t(q,p)=(q,p-t\cdot f'(q)),
\end{equation}
where $q$ is a coordinate on $\mathbb{S}^1$ and $p$ is a coordinate on a fibre. We denote by $H$ the autonomous Hamiltonian obtained by cutting of $H_f$ outside some compact set so that $H_f\equiv H$ on a compact neighbourhood $K$  of $\phi_{H_f}^1(o_{\mathbb{S}^1})$.  Its flow on $K$, denoted by $\phi^t$, is given in~(\ref{eq:Ham1}). We see that the $n$-th power of $\phi$ is given by
$$
\phi^n(q,p)=(\phi^1)^n(q,p)=(q,p-n\cdot f'(q)),
$$
and it is a time-one map generated by $H^{\sharp n}(q,p)=nf(q)$ on the compact $K$. 
A generator of $HF_*(o_{\mathbb{S}^1}, T^*_{x_1}\mathbb{S}^1:H^{\sharp n})$ is the path
$$
\gamma(t) = \phi^{nt}\left((\phi^n)^{-1}(z)\right) = 
\phi^{nt}(x_1,0) = (x_1,-ntf'(x_1)) = (x_1,0),
$$
where $z=\phi^n(x_1,0)\in\phi^n(o_{\mathbb{S}^1})$. Since $\ell_+^{\{x_1\}}(\phi^n)$ is a critical value of the action functional we conclude
$$\ell_+^{\{x_1\}}(\phi^n)=\mathcal{A}_{H^{\sharp n}}(\gamma)=-\int_\gamma p\,dq+\int_0^1H^{\sharp n}(\gamma(t),t)\,dt=n f(x_1)
$$
and 
$$
\sigma^{\{x_1\}}(\phi)=f(x_1)=\frac{1}{n}\ell^{\{x_1\}}_+(\phi^n).
$$
For the same reason, for the other critical point $x_2$ we have 
$$
\sigma^{\{x_2\}}(\phi)=f(x_2)>f(x_1)=\sigma^{\{x_1\}}(\phi).
$$

Let us briefly discuss the proof in~\cite{MVZ} of the fact that the quasi-morphism constructed therein coincides with Viterbo's homogenization on a torus. If we reformulate it in terms of the action functional we can see that the key step is the equality
\begin{equation}\label{MVZJedn}
\frac{1}{n}\ell_+(\phi_H^n)=\ell_+(\phi_{H_n}^1)
\end{equation}
that holds for every $n$. Here, $\ell_+:=\ell_+^M$.  The Hamiltonian $H_n$ is a homogenization of the Hamiltonian $H$ on a torus and is given by $H_n(q,p)=H(nq,p)$. 

The following natural question is arrising here. Is it true that the relation~(\ref{MVZJedn}) holds for any conormal spectral invariants $\ell_+^N$, i.e. does it hold
\begin{equation}\label{Quest}
\frac{1}{n}\ell_+^N(\phi_H^n)=\ell_+^N(\phi_{H_n}^1)
\end{equation}
for a closed submanifold $N$? If we consider the situation we described above, when $M$ is a circle, the Hamiltonian $H$ is the lift of the Morse function $f$ and $N=\{x_1\}$ is the minimum point of $f$ we can see that the left-hand side in~(\ref{Quest}) is equal to $f(x_1)$. We compute the right-hand side  as follows. The Hamiltonian $H_n$ generates the Hamiltonian flow
$$
\phi_{H_n}^t(q,p)=(q,p-tnf'(nq)).
$$
A generator of $HF_*(o_{\mathbb{S}^1}, T^*_{x_1}\mathbb{S}^1:H_n)$ is
$$
v(t)=\phi_{H_n}^t((\phi_{H_n}^1)^{-1}(z))
$$
where $z=\phi_{H_n}^1(x_1,0)$. Then 
$$v(t)=(x_1,-tnf'(nx_1)). 
$$
The corresponding critical value of the action functional is
$$\ell_+^{\{x_1\}}(\phi_{H_n}^1)=\mathcal{A}_{H_n}(v)=-\int_v p\,dq+\int_0^1H_n(v(t),t)\,dt=f(nx_1).
$$
It follows that the relation~(\ref{Quest}) does not hold. Moreover, the $\overline{\lim}$ of its two sides do not necessarily coincide. For example, if in the representation of $\mathbb{S}^1$ as the quotient $[0,1]/\{0,1\}$ the minimum point $x_1$ is an irrational number, then the orbit $\{nx_1\}$ is dense in $\mathbb{S}^1$. Therefore, in that case
$$
\overline{\lim}\, \ell^{\{x_1\}}_+(\phi^1_{H_n}) =
\lim_{m\to 0}\sup_{n\geq m} f(nx_1) = \max f,
$$
while 
$$
\frac{1}{n}\ell^{\{x_1\}}_+(\phi^N_H)=f(x_1)=\min f.
$$

\section{Proof of the main theorems}\label{sec:proof}

In this section we prove Theorem~\ref{ThrmQuasiMorph} and Theorem~\ref{ThrmQuasiState}. We keep the notations introduced in previous sections and all of the properties from therein. 

\subsection{Proof of Theorem~\ref{ThrmQuasiMorph}} Since the sequence 
$\frac{a_n}{n}$ is bounded,  $\overline{\lim}\frac{a_n}{n}\in\mathbb{R}$, and thus 
$\overline{\lim}\frac{\ell_+^N(\phi^n)}{n}\in\mathbb{R}$.  It means that 
$\sigma^N$ is a well defined map. 

First we prove the property~(\ref{Homog}). We easily obtain one inequality
\begin{equation}\label{1side}
\begin{aligned} 
\sigma^N(\phi^l)&=\overline{\lim\limits_n}\frac{\ell_+^N(\phi^{nl})}{n}\\
&=l\cdot\overline{\lim\limits_n}\frac{\ell_+^N(\phi^{nl})}{nl}\\
&\le l\cdot\overline{\lim\limits_n}\frac{\ell_+^N(\phi^{n})}{n}=l\,\sigma^N(\phi).
\end{aligned}
\end{equation}
The opposite inequality follows by the following observation. Let $\{a_{n_k}\}$ be a subsequence of $a_n$ such that 
$$\overline{\lim\limits}\frac{a_n}{n}=\lim\limits_{k\to+\infty}\frac{a_{n_k}}{n_k}.
$$
We define a sequence $\{m_k\}$ that is the quotient when dividing $n_k$ by $l$
$$n_k=l\cdot m_k+r_k,
$$
and $r_k\in\{0,1,\ldots,l-1\}$ is the remainder.  It is obvious that $m_k\to+\infty$ when $n_k\to+\infty$. Using the fact that $a_n$ is non-decreasing sequence we get
$$
\begin{aligned}
\frac{a_{n_k}}{n_k}&=\frac{a_{l\cdot m_k+r_k}}{l(m_k+1)}\frac{l(m_k+1)}{n_k}\\
&\le\frac{a_{l\cdot m_k+l}}{l(m_k+1)}\frac{l(m_k+1)}{n_k}\\
&=\frac{a_{l(m_k+1)}}{l(m_k+1)}\frac{l(m_k+1)}{n_k}.
\end{aligned}
$$
Passing to the $\overline\lim$ we get
$$
\begin{aligned}
\overline{\lim\limits_k}\frac{a_{n_k}}{n_k}&\le \overline{\lim\limits_k}\left[\frac{a_{l(m_k+1)}}{l(m_k+1)}\frac{l(m_k+1)}{n_k}\right]\\
&=\overline{\lim\limits_k}\frac{a_{l(m_k+1)}}{l(m_k+1)},
\end{aligned}$$
where the last equality holds because $\frac{l(m_k+1)}{n_k}=\frac{n_k-r_k+l}{n_k}\to1$ when $k\to+\infty$.
Since $l(m_k+1)$ is a subsequence of $nl$ it holds
\begin{equation}\label{help}
\begin{aligned}
\overline{\lim\limits_n}\frac{a_{nl}}{nl}&\ge\overline{\lim_k}\frac{a_{l(m_k+1)}}{l(m_k+1)}\\
&\ge\overline{\lim\limits_k}\frac{a_{n_k}}{n_k}=\lim\limits_{k\to+\infty}\frac{a_{n_k}}{n_k}.
\end{aligned}
\end{equation}
We are now back to our map $\sigma^N$. Using the  subsequence $n_k$  defined 
above we see that
$$
\begin{aligned}
\sigma^N(\phi)&=\overline{\lim}\frac{\ell^N_+(\phi^n)}{n}\\
&=\overline{\lim}\frac{a_n}{n}-\ell^M_+(\phi^{-1})\\
&=\lim\limits_{k\to+\infty}\frac{a_{n_k}}{n_k}-\ell^M_+(\phi^{-1}).
\end{aligned}
$$
Next step is to correlate $a_{nl}$ and $\sigma^N(\phi^{l})$
\begin{equation}\label{2side}
\begin{aligned}
\sigma^N(\phi^l)&=\overline{\lim\limits_n}\frac{\ell_+^N(\phi^{nl})}{n}\\
&=\overline{\lim\limits_n}\frac{a_{nl}-nl\cdot \ell^M_+(\phi^{-1})}{n}\\
&=\overline{\lim\limits_n}\frac{a_{nl}}{n}-l\cdot \ell^M_+(\phi^{-1})\\
&=l\cdot\overline{\lim\limits_n}\frac{a_{nl}}{nl}-l\cdot \ell^M_+(\phi^{-1})\\
&\stackrel{(\ref{help})}\ge l\cdot\lim\limits_{k\to+\infty}\frac{a_{n_k}}{n_k}-l\cdot \ell^M_+(\phi^{-1})\\
&=l\cdot(\sigma^N(\phi)+\ell^M_+(\phi^{-1}))-l\cdot \ell^M_+(\phi^{-1})\\
&=l\cdot\sigma^N(\phi).
\end{aligned}\end{equation}

The property~(\ref{Homog}) follows from~(\ref{1side}) and~(\ref{2side}).

Now we prove the property~(\ref{ConjInvar}).  The inequality~(\ref{triang1}) implies that for all $\psi,\phi\in\Ham(T^*M)$ we have
$$
\begin{aligned}
\ell_+^N(\psi\phi^n\psi^{-1})&\le \ell_+^M(\psi)+\ell_+^N(\phi^n\psi^{-1})\\
&\le \ell^M_+(\psi)+\ell^M_+(\psi^{-1})+\ell_+^N(\phi^n),
\end{aligned}
$$
and also
$$\begin{aligned}
\ell_+^N(\phi^n)&=\ell_+^N(\psi^{-1}\psi\phi^n\psi^{-1}\psi)\\
&\le \ell^M_+(\psi^{-1})+\ell^M_+(\psi)+\ell_+^N(\psi\phi^n\psi^{-1}).
\end{aligned}
$$
We bound the difference $\ell_+^N(\psi\phi^n\psi^{-1})-\ell_+^N(\phi^n)$ 
from both sides
$$-\ell^M_+(\psi)-\ell_+^M(\psi^{-1}) \le 
\ell_+^N(\psi\phi^n\psi^{-1})-\ell_+^N(\phi^n)\le
\ell^M_+(\psi)+\ell_+^M(\psi^{-1}). 
$$
Dividing these three inequalities by $n$ we have
$$
0=\lim\limits_{n\to+\infty} \left[\frac{\ell_+^N(\psi\phi^n\psi^{-1})}{n}-
\frac{\ell_+^N(\phi^n)}{n}\right].
$$
Using the property of $\overline{\lim}$ and $\underline{\lim}$ we conclude that
$$\begin{aligned}
0&=\overline{\lim}\left[\frac{\ell_+^N(\psi\phi^n\psi^{-1})}{n}-\frac{\ell_+^N(\phi^n)}{n}\right]\\
&\ge\overline{\lim}\frac{\ell_+^N(\psi\phi^n\psi^{-1})}{n}+
\underline{\lim}\left(-\frac{\ell_+^N(\phi^n)}{n}\right)\\
&=\overline{\lim}\frac{\ell_+^N(\psi\phi^n\psi^{-1})}{n}-
\overline{\lim}\frac{\ell_+^N(\phi^n)}{n}
\end{aligned}
$$
and
$$\begin{aligned}
0&=\overline{\lim}\left[\frac{\ell_+^N(\phi^n)}{n}-\frac{\ell_+^N(\psi\phi^n\psi^{-1})}{n}\right]\\
&\ge \overline{\lim}\frac{\ell_+^N(\phi^n)}{n}+
\underline{\lim}\left(-\frac{\ell_+^N(\psi\phi^n\psi^{-1})}{n}\right)\\
&=\overline{\lim}\frac{\ell_+^N(\phi^n)}{n}-
\overline{\lim}\frac{\ell_+^N(\psi\phi^n\psi^{-1})}{n}.
\end{aligned}$$
Thus
$$
\overline{\lim}\frac{\ell_+^N(\phi^n)}{n}=
\overline{\lim}\frac{\ell_+^N(\psi\phi^n\psi^{-1})}{n}.
$$
The conjugation-invariance of $\sigma^N$ now easily follows since
$$
\begin{aligned}
\sigma^N(\psi\phi\psi^{-1})&=
\overline{\lim}\frac{\ell_+^N((\psi\phi\psi^{-1})^n)}{n}\\
&=\overline{\lim}\frac{\ell_+^N(\psi\phi^n\psi^{-1})}{n}\\
&=\overline{\lim}\frac{\ell_+^N(\phi^n)}{n}=\sigma^N(\phi),
\end{aligned}
$$
for every $\phi,\psi\in\Ham(T^*M)$.

The proof of the point~(\ref{Lipsic}) is the same as the proof of the property (iii) in Theorem 1.3 in~\cite{MVZ}, since conormal spectral invariants $\ell_+^N$ satisfy the same inequality~(\ref{LipsicInvar}) that we have already mentioned.

Arguing in a similar way as the authors in~\cite{MVZ} we obtain the claim~(\ref{Razdvojivo}). It follows from Proposition 2.17 in~\cite{MVZ}  that for every $\psi\in\Ham(T^*M)$ that displaces an open subset $U\subset T^*M$ and every $\phi\in\Ham_U(T^*M)$ it holds
$$-\Gamma(\psi)\le \ell^M_+(\phi)\le\Gamma(\psi).
$$
Here $\Gamma(\psi)$ is the spectral norm of $\psi$, defined using Hamiltonian spectral invariants $c_{\pm}$ by $\Gamma(\psi):=c_+(\psi)-c_-(\psi)$
(see Section 2.2 in~\cite{MVZ} for more details). We already know that
$$
\ell_+^N(\phi)\le \ell^M_+(\phi)
$$
and  also
$$
0=\ell_+^N(\Id)=\ell_+^N(\phi^{-1}\phi)\le \ell_+^M(\phi^{-1})+\ell_+^N(\phi).
$$
Since $\phi^{-1}$ has a compact support in $U$ combining previous inequality we get
$$
-\Gamma(\psi)\le-\ell_+^M(\phi^{-1})\le \ell_+^N(\phi)\le\Gamma(\psi).
$$
Every iteration $\phi^n$ belongs to $\Ham_U(T^*M)$ thus
$$\overline{\lim}\frac{\ell_+^N(\phi^n)}{n}=
\lim\limits_{n\to+\infty}\frac{\ell_+^N(\phi^n)}{n}=\lim\limits_{n\to+\infty}\frac{\Gamma(\psi)}{n}=0.
$$

Now we prove the property~(\ref{Fragmen}). Let us suppose that $\phi=\phi_H^1$ where $supp(H)\subset U$ and $U$ is the element of $\mathcal{U}$. If $\psi$ is any element of $\Ham(T^*M)$ then $\phi_j:=\psi^j\phi\psi^{-j}$ is dominated by some element in $\mathcal{U}$ (see~\cite{MVZ}, p. 235 and p. 207 for the definition). Since $(\phi\psi)^n=\phi_0\phi_1...\phi_{n-1}\psi^n$ using the triangle inequality we conclude
$$\begin{aligned}
\ell^N_+((\phi\psi)^n)&=\ell^N_+(\phi_0\phi_1...\phi_{n-1}\psi^n)\\
&\le\ell^M_+(\phi_0\phi_1...\phi_{n-1})+\ell^N_+(\psi^n)\\
&\le\sum\limits_{j=0}^n \ell^M_+(\phi_j)+\ell^N_+(\psi^n)\\
&\le ne(\mathcal{U})+\ell^N_+(\psi^n).
\end{aligned}
$$
The last inequality follows from the fact that $\ell^M_+(\phi_j)\le e(\mathcal{U})$ for every $j$ since $\phi_j$ is generated by Hamiltonian whose support is a displaceable subset with displacement energy smaller then $e(\mathcal{U})$. Taking the limsup of the both side in inequality
$$
\frac{\ell^N_+((\phi\psi)^n)}{n}\le e(\mathcal{U})+\frac{\ell^N_+(\psi^n)}{n}
$$
we get the estimate 
$$\sigma^N(\phi\psi)\le\sigma^N(\psi)+e(\mathcal{U}).
$$
The point~(\ref{Fragmen}) now follows by induction on the positve integer $\|\phi\|_{\mathcal{U}}$.

The point~(\ref{Commute}) follows again from the triangle inequality
$$
\ell_+^N((\phi\psi)^n)=\ell_+^N(\phi^n\psi^n)\le 
\ell_+^N(\phi^n)+\ell_+^M(\psi^n),
$$
dividing by $n$ and taking the upper limit of the both sides.
\qed

\subsection{Proof of Theorem~\ref{ThrmQuasiState}}
The claims~(\ref{NormState}) -~(\ref{VanState}) follow directly from Theorem~\ref{ThrmQuasiMorph}. The point~(\ref{QAddState}) follows from the point~(\ref{Commute}) of Theorem~\ref{ThrmQuasiMorph}:
$$
\zeta^N(H+K)\le\zeta^N(H)+\zeta^M(K)=\zeta^N(H),
$$
since $\zeta^M(K)=0$ when $\supp(K)$ is displaceable (property~(\ref{Razdvojivo}) in Theorem~\ref{ThrmQuasiMorph}). Since
$\{H+K,-K\}=0$ and the support of $-K$ is displaceable as well  we obtain 
$$\zeta^N(H)=\zeta^N(H+K+(-K))\le\zeta^N(H+K).
$$\qed

\section{A comment on manifolds with boundary}\label{sec:boundary}

In the proof of the previous theorems, we have shown that it is possible to construct quasi-morphisms using Lagrangian spectral numbers even in the case when we do not have a product in homology that ensures that they are subadditive. With the help of constructions in~\cite{KO,KMN1,KMN2}, we can apply these conclusions to the case when N is a manifold with a boundary.

Let $M$ be a closed  manifold as above, and $N\subset M$ a submanifold with boundary $\partial N$.  There are two possible choices of a definition of the singular Lagrangian submanifold $\overline{\nu}^*N$:
$$
\overline{\nu}^*_{\pm}N:=\nu^*_{\pm}(\partial N)\cup\nu^*(\mathrm{Int\,N}),
$$
where
$$
\nu^*_{\pm}(\partial N):=\{(q,p)\in\nu^*N \mid \pm p(\vec{n})\leq 0\}
$$
for $\vec{n}\in T|_{\partial N}N$ inner normal to $\partial N$. Lagrangian Floer homology with boundary conditions $(o_M,\overline{\nu}^*_{\pm}N)$ was first constructed by Kasturirangan and Oh~\cite{KO} for open subset $N$ (see~\cite{KMN2} for the generalization to submanifolds of positive codimension). The two choices above give rise to two Lagrangian Floer homologies, isomorphic to singular homologies $H_*(N)$ and $H_*(N,\partial N)$.  For a homology class $\alpha\in H_*(N)$ (or a relative homology class $\alpha\in H_*(N,\partial N)$, there are well defined spectral numbers $\ell^N_+(\alpha,\cdot)$  defined as a minimax values in corresponding Lagrangian Floer homology (see~\cite{KMN1,KMN2} for details). 

If $\dim\,N=\dim\,M$ then there exists a pair-of-pants product on Floer homology that intertwines with the intersection product in the homology and the PSS isomorphism (see Theorem 3 in~\cite{KMN1}).  As a consequence, one can define the subaditive spectral invariants $\ell^N_+(\phi)$ by using the idempotent element $[N]\in H_{\dim\,N}(N,\partial N)$ and setting 
$\ell^N_+(\cdot):=\ell^N_+([N],\cdot)$ . If $\dim\,N<\dim\,M$ such a product does not exist, but there is a well defined exterior product $\widetilde{\cap}$ analogous to~(\ref{eq:extcap}) (see Theorem B in~\cite{KMN2}). Arguing as before we can conclude that in both cases there are well defined partial quasi-morphism
$$
\sigma^N_{\partial}:\mathrm{Ham}(T^*M)\rightarrow\mathbb{R},\quad \sigma^N_{\partial}(\phi):=\overline{\lim}\frac{\ell^N_+(\phi^n)}{n}
$$
(if $\dim\,N=\dim,M$ the upper limit above is just the limit), 
and a partial quasi-state
$$
\zeta^N_{\partial}:C^{\infty}_c(T^*M)\rightarrow\mathbb{R}, \quad \zeta^N_{\partial}(H):=\sigma^N_{\partial}(\phi^H_1).
$$

\end{document}